\documentclass{siamltex}[final]
\usepackage{amsfonts,amsmath,color,verbatim}
\usepackage{stmaryrd}
\usepackage[mathscr]{eucal}
\usepackage{epsfig}
\usepackage{comment}
\usepackage[ruled,vlined]{algorithm2e}
\usepackage{color}
\usepackage{url}
\usepackage{GGGraphics}

\def\d{{\mathrm d}}
\def\e{{\mathrm e}}

\def\eps{\varepsilon}

\def\calM{{\mathcal M}}
\def\calJ{{\mathcal J}}

\newdimen\GGGlength
\newdimen\GGGheight
\newbox\GGGbox
\def\GGGput[#1,#2](#3,#4)#5{%
  \setbox\GGGbox\vbox{\hbox{#5}\kern0pt}%
  \GGGlength\wd\GGGbox%
  \divide\GGGlength by100 \multiply\GGGlength by#1%
  \GGGheight\ht\GGGbox%
  \divide\GGGheight by100 \multiply\GGGheight by#2%
  \put(#3,#4){\kern-\GGGlength\raise-\GGGheight\box\GGGbox}}

\newcommand{\R}{\mathbb R}

\newcommand{\bq}{\begin{equation}}
\newcommand{\eq}{\end{equation}}

\def\e{{\rm e}}
\def\d{{\rm d}}
\def\bigo{{\mathcal O}}

\def\diag{{\mathrm {Diag}}\,}

\def\1{{\bf 1}}

\title{Applying stiff integrators for ODEs and DDEs to problems with distributed delays}

\author{Nicola Guglielmi\footnotemark[1] \and Ernst Hairer\footnotemark[2]}
\begin{document}

\maketitle

\renewcommand{\thefootnote}{\fnsymbol{footnote}}
\footnotetext[1]{Division of Mathematics, 
Gran Sasso Science Institute,
Via Crispi 7,
I-67100   L'Aquila,  Italy. Email: {\tt nicola.guglielmi@gssi.it}}
\footnotetext[2]{Section de Math{\'e}matiques, Universit\'e de Gen{\`e}ve,
       CH-1211 Gen{\`e}ve 24, Switzerland.\\
       Email: {\tt ernst.hairer@unige.ch}}
\renewcommand{\thefootnote}{\arabic{footnote}}

%
%
%
%
%
%
%
%
%
%
%
%
\begin{abstract}
There exist excellent codes for an efficient numerical treatment of stiff and
differential-algebraic problems. Let us mention {\sc Radau5} which is based on the
$3$-stage Radau IIA collocation method, and its extension to problems with
discrete delays {\sc Radar5}. The aim of the present work is to present a technique
that permits a direct application of these codes to problems having a
right-hand side with an additional
distributed delay term (which is a special case of an integro-differential equation).
Models with distributed delays are of increasing importance
in pharmacodynamics and pharmacokinetics for the study of the interaction between drugs and the
body.

The main idea is to approximate the distribution kernel of the integral term by
a sum of exponential functions or by a
quasi-polynomial expansion, and then to transform the distributed (integral) delay
term into a set of ordinary differential equations. This set
is typically stiff and, for some distribution kernels (e.g., Pareto
distribution), it contains discrete delay terms with constant delay.
The original equations augmented by this set of ordinary differential equations
can have a very large dimension, and a careful treatment of the solution of the
arising linear systems is necessary.

The use of the codes {\sc Radau5} and {\sc Radar5} is illustrated at three examples
(two test equations and one problem taken from pharmacodynamics). The driver programs
for these examples are publicly available from the homepages of the authors.

\end{abstract}



\pagestyle{myheadings} \thispagestyle{plain}
\markboth{N. GUGLIELMI,  E. HAIRER}{Applying codes for ODEs and DDEs to problems with distributed delays}


\section{Introduction}
\label{sect:intro}

The motivation of the present work is in dealing numerically with mathematical models 
consisting of delay integro-differential equations, as those arising in 
pharmacokinetics and pharmacodynamics, where
the interaction between drugs and the body is studied (e.g.,
\cite{wesolowski20cot,Krzyzanski18,Humphries18}). 
The problem is usually modeled by ordinary
or delay, differential or differential-algebraic equations, which typically contain
integral terms of convolution type (also termed distributed delay).
Differential equations with a distributed delay term are a special case of
Volterra integro-differential equations. There is a large literature on the numerical
discretization of such problems, and their accuracy and stability is well investigated.
Let us mention the monographs by H.\ Brunner
\cite{brunner04cmf,brunner86tns}.
A different approach for problems with weakly singular kernel is the use of
discretized fractional calculus by C.\ Lubich \cite{lubich86dfc} (see also \cite{hairer87nmf}),
and  the oblivious convolution quadrature \cite{Sc06,LF08}
for more general kernels.
Available codes mostly make use of constant step size.

Implementing an integration method
for the considered class of problems is more or less straight-forward,
if constant time steps are used.
However, in realistic situations one is often confronted
with initial layers (like the problem
of Section~\ref{sect:chemo}) and/or fast transitions between different states, so that
flexibility in the choice of step size is an essential ingredient
for efficiency.
Although there exist excellent codes for solving nonstiff and stiff ordinary differential equations,
codes for differential-algebraic equations, and codes for delay differential equations
equipped with sophisticated step size strategies,
we are not aware of a code based on a variable step size
time integrator that can efficiently treat problems with
distributed delays.

Instead of extending a numerical integrator for Volterra integro-differential equations
to a variable step size code, we propose to change the problem in such a way
that any code for stiff and differential-algebraic (delay) equations
can be applied.
The idea is to approximate the distribution kernel of the integral term by
a sum of exponential functions (multiplied by a polynomial) 
and then to transform the integral delay
term into a set of differential equations. The use of an approximation by a sum of exponential
functions is not new in numerical analysis. In \cite{braess09ote} it is applied
for an efficient computation of high-dimensional integrals.
In \cite{LF08} it is used for developing a variable step size integrator that is applied
to interesting problems like a blow-up problem for a nonlinear Abel integral equation
and a fractional diffusion-reaction system.
New in the present work is that we do not propose another time integrator, but we show
that a large class of problems with distributed delay can be solved by standard
software for stiff and differential-algebraic (delay) equations.
This is illustrated at the hand of {\sc Radau5} \cite{hairer96sod}
for stiff and differential-algebraic equations, and its extension
{\sc Radar5} \cite{guglielmi01irm,guglielmi08cbp} for problems including delay arguments.
New is also an 
application to problems in
pharmacokinetics and pharmacodynamics.
\medskip

\noindent
{\it Outline of the paper.}
Section \ref{sect:dgln} introduces the class of problems that can be treated with the
codes presented in this work. It comprises nonstiff, stiff, differential-algebraic,
and delay differential equations, with the additional feature that distributed delay
terms are included.
In approximating the kernel of a distributed delay term by a sum of exponential
functions multiplied by a polvnomial,
the problem is transformed into a system
without distributed delay terms.
Since the resulting system is typically of a much larger dimension than the original problem,
Section \ref{sect:radar} presents an algorithm
that reduces considerably the complexity of the required solution of linear systems.
Furthermore, the stability of the proposed algorithm and the determination
of the accuracy parameters are discussed.
Section \ref{sect:sumofexp} describes the
approach by Beylkin and Monz\' on \cite{beylkin10abe} for approximating the factor $t^{-\alpha}$
(typically present in weakly singular kernels) by a sum of exponential functions.
Parameters in the approximation are selected to keep the approximation error under a
given level.
Section \ref{sect:distributed} presents important examples of distributed delays
(the gamma distribution and the Pareto distribution) commonly used in pharmacodynamics and pharmacokinetics.
Finally, Section~\ref{sect:numerical} provides numerical evidence of the efficiency of the
proposed approach by applying the codes {\sc Radau5} and {\sc Radar5}
to two test examples and to an example from chemotherapy-induced myelosuppression.
The codes together with drivers for problems with distributed delays are made publicly available.

\section{Differential equations with distributed delay}
\label{sect:dgln}

We consider differential equations of the form
\begin{equation}\label{problem}
M \dot y (t) = f\bigl( t, y(t), y(t - \tau ) , I (y)(t)\bigr),
\qquad y(0)=y_0 ,\quad y(t) = \eta (t) ~~ \text{for} ~~ t<0 ,
\end{equation}
where $y \in \R^d$, $M$ is a constant $d\times d$ matrix, the delay satisfies $\tau \ge 0$, and
\begin{equation}\label{distridelay}
I (y)(t) = \int_0^t k(t-s)g\bigl( s, y(s)\bigr) \,\d s 
\end{equation}
is a distributed delay term.
The vector functions $\eta(t)$, $f(t,y,v,I)$, and the scalar functions $k(t)$, $g(t,y)$
are assumed to be smooth. The kernel $k(t)$ is allowed to have an integrable singularity
at the origin.

This formulation contains ordinary differential equations,
differential-algebraic equations (for singular $M$), delay differential equations, and
integro-differential equations as special cases. All considerations of the present work
extend straight-forwardly to the case, where the delay $\tau = \tau (t, y(t))\ge 0$ is
time and state dependent, and where
several delay terms
and several integral terms $I_i(t)(y)$ (with different kernels $k_i(t)$ and different
functions $g_i(t,y)$) are present. For notational convenience we restrict our study
to the formulation \eqref{problem}.
We focus on the situation, where at least one integral term
\eqref{distridelay} is present in \eqref{problem}.

Assume that we have at our disposal a reliable code for solving differential equations
of the form \eqref{problem}, where no integral term is present. Let us mention the
code {\sc Radau5} of \cite{hairer96sod} for stiff and differential-algebraic
equations and its extension {\sc Radar5} of \cite{guglielmi01irm} for delay
differential-algebraic equations. Our aim is to demonstrate how such codes can be applied
to the solution of \eqref{problem} with integral terms included, without changing them.
The main idea is to approximate the kernel $k(t)$ by a finite sum of exponential
functions multiplied by a polynomial. We thus assume that
\begin{equation}\label{sum-exp-pol}
k(t) = \sum_{i=1}^n p_i(t) \,\e^{-\gamma_i t},\qquad  p_i(t)= \sum_{j=0}^{m_i} c_{i,j}\, t^j
\end{equation}
with real coefficients $c_{i,j}$ and $\gamma_i$, so that the distributed delay term
\eqref{distridelay} can be written as
\begin{equation}\label{def-I}
I (y)(t) = \sum_{i=1}^n \sum_{j=0}^{m_i} c_{i,j}  z_{i,j} (t), \qquad z_{i,j} (t) =
\int_0^t (t-s)^j \,\e^{-\gamma_i (t-s)} g\bigl( s, y(s)\bigr) \,\d s .
\end{equation}
Differentiation of $z_{i,j}(t)$ with respect to time yields, for $i=1,\ldots ,n$,
\begin{equation}\label{zij-dgl}
\dot z_{i,j} (t) = - \,\gamma_{i}\, z_{i,j}(t) 
+ \Bigg\{\begin{array}{cl} g\bigl( t, y(t)\bigr)  & j=0 \\[2mm]
j \, z_{i,j-1}(t)  & j=1,\ldots ,m_i .\end{array}
\end{equation}
We insert \eqref{def-I} into the differential equation \eqref{problem}, we denote
$z_i = \bigl(z_{i,0},\ldots ,z_{i,m_i}\bigr)^\top$, and we add
these differential equations for $z_{i,j}(t)$ to the original equations.
This gives the augmented system of differential equations
\begin{equation}\label{problem-aug}
\hskip -0.1cm
\begin{array}{rcll} M \dot y (t) & \!\!\!=\!\!\! & \displaystyle
f\Bigl( t, y(t), y(t - \tau ) , 
 \sum_{i=1}^n \sum_{j=0}^{m_i} c_{i,j}  z_{i,j} (t)\Bigr) ,&
\ \begin{array}{l} y(0)=y_0  \\[1mm] y(t) =  \eta (t) ~~ \text{for} ~~ t<0 \end{array}
\\[6mm]
\dot z_1(t) & \!\!\!=\!\!\! & J_1 z_1(t) + e_1 g\bigl( t , y(t)\bigr) & \ z_1(0) = 0 \\[1mm]
& \!\!\!\! \cdots \!\!\!\!&&  \qquad \cdots \\[1mm]
\dot z_n(t) & \!\!\!=\!\!\! & J_n z_n(t) + e_n g\bigl( t , y(t)\bigr) & \ z_n(0) = 0 
\end{array}
\end{equation}
for the super vector $\bigl(y,z_1,\ldots , z_n\bigr)$ of dimension $d+(m_1+1)+\ldots + (m_n+1)$.
Here, we use the notation
\[
J_i = \begin{pmatrix}
-\gamma_i  & 0 & \cdots &\cdots &  0 \\[-1.5mm]
1 & -\gamma_i &\ddots & &\vdots \\[-1.5mm]
0 & 2 & -\gamma_i & \ddots &  \vdots \\[-1.5mm]
\vdots & \ddots & \ddots & \ddots & 0 \\[-0.5mm]
0 &\cdots & 0 & m_i & -\gamma_i \,
\end{pmatrix} \in \R^{m_i+1,m_i+1},\qquad e_i =  \begin{pmatrix}
1\\
0 \\[-1.5mm]
  \vdots \\[-1.5mm]
\vdots\\[-0.5mm]
0
\end{pmatrix} \in \R^{m_i+1}
\]
for the 
lower bidiagonal matrix and the first unit vector, both of dimension $m_i+1$.
This
gives a formulation, where also those codes can be applied, that are not adapted
to the treatment of distributed delays \eqref{problem}.

\section{Solving distributed delay equations via the augmented system}
\label{sect:radar}

For the case that the original system \eqref{problem} is nonstiff and the exponents
$\gamma_i$ in \eqref{sum-exp-pol} are of moderate size, any code for nonstiff ODEs
or DDEs can be applied to the augmented system \eqref{problem-aug}.
However, in important applications (see Section~\ref{sect:sumofexp}) some of the
$\gamma_i$ can be very large, so that the augmented system is stiff independent of
the behaviour of \eqref{problem}. Therefore, a stiff solver (e.g., {\sc Radau5} for the
case of ODEs, and {\sc Radar5} in the presence of retarded arguments) has to be applied
to \eqref{problem-aug}. Stiff solvers are in general implicit and require an efficient
solution of linear systems of the form (see \cite[Section IV.8]{hairer96sod})
\begin{equation}\label{lin-syst}
\Bigl( (\gamma h )^{-1} \calM  - \calJ \Bigr) u = a ,
\end{equation}
where $a = (a_0, a_1,\ldots , a_n)^\top$ is a given vector and $u =(u_0,u_1,\ldots , u_n)^\top$
is the solution of the system. Here,~$\gamma$ is a real or complex parameter of the stiff integrator,
$h$ is the time step size,
$\calM = \diag \bigl( M, I_1, \cdots , I_n\bigr)$ is a block diagonal matrix
with $I_i$ being the identity matrix of dimension $m_i+1$, and
the Jacobian matrix of the augmented system is
\begin{equation}\label{jacobian}
\calJ = \begin{pmatrix}
J  & C_1 & C_2 &\cdots &  C_n \\[1.mm]
B_1 & J_1 & 0 & \cdots & 0 \\[-0.5mm]
B_2 & 0 & J_2 & \ddots & \vdots \\[-0.5mm]
\vdots & \vdots & \ddots & \ddots & 0 \\[0.5mm]
B_n & 0 & \cdots & 0 & J_n
\end{pmatrix}  ,
\end{equation}
where $J = (\partial f/\partial y)(t,y,v,I)$ is a square matrix of dimension $d$,
and the rank-one matrices $B_i$ and $C_i$ are given as follows:
$B_i = e_i  g_y^\top$ with $g_y = (\partial g/\partial y)(t,y)$, and
$C_i =f_I c_i^\top$ with $f_I = (\partial f/\partial I)(t,y,v,I)$ and
$c_i = (c_{i,0},c_{i,1},\ldots ,c_{i,m_i})^\top$. All derivatives are
evaluated at the current integration point.

\subsection{Efficient solution of the linear system \eqref{lin-syst}}
\label{sect:lin-syst}

The linear system \eqref{lin-syst} is of dimension $d+N$, where $N = (m_1+1) +
\ldots + (m_n+1)$
can be much larger than $d$, the dimension of the original problem. Without exploiting
the special structure of \eqref{lin-syst} its solution costs $\bigo \bigl((d+N)^3\bigr)$
flops. We propose to solve the linear system in the following way.

The first block of \eqref{lin-syst} is equivalent to
\begin{equation}\label{lin-syst1}
\Bigl( (\gamma h)^{-1} M - J \Bigr) u_0 = a_0 + \sum_{i=1}^n C_i u_i ,
\end{equation}
and the other blocks give 
\begin{equation}\label{lin-systi}
\bigl( (\gamma h)^{-1} I_i - J_i \bigr) u_i = a_i + B_i u_0 .
\end{equation}
This equation permits to express $u_i$ in terms of $u_0$. Inserted into \eqref{lin-syst1}
yields
\begin{equation}\label{lin-syst1a}
\Bigl( (\gamma h)^{-1} M - J \Bigr) u_0 = a_0 + \sum_{i=1}^n C_i 
\Bigl( (\gamma h)^{-1} I_i - J_i \Bigr)^{-1} \bigl( a_i + B_i u_0 \bigr) ,
\end{equation}
which can be written as
\begin{eqnarray} \nonumber
&& \Bigl( (\gamma h)^{-1} M - \widehat J \, \Bigr) u_0 = a_0 + f_I\sum_{i=1}^n c_i^\top
\Bigl( (\gamma h)^{-1} I_i - J_i \Bigr)^{-1}\!  a_i , \\[2mm]
&& \widehat J = J + f_I g_y^\top \sum_{i=1}^n c_i^\top
\Bigl( (\gamma h)^{-1} I_i - J_i \Bigr)^{-1}\! e_i .
\label{lin-syst1b}
\end{eqnarray}
Note that $\widehat J$ is a rank-one perturbation of $J$. Since
$\bigl( (\gamma h)^{-1} I_i - J_i \bigr)$ is a lower bidiagonal matrix of dimension~$m_i+1$,
the computation of the constant in the rank-one perturbation and the
computation of the right-hand side of
the linear system \eqref{lin-syst1b} require not more than $\bigo (N)$ flops.
Having computed $u_0$, the $u_i$ are obtained from \eqref{lin-systi}.
All in all, this gives an algorithm that has a cost of
$\bigo (d^3 ) + \bigo (N) $ flops.

\subsection{Stable approximation of the distributed delay term: exponential case}
\label{sect:stable}

It may happen that several $\gamma_n$ in the kernel approximation \eqref{sum-exp-pol}
are very large and positive, so that the summands are not well scaled. Consider,
for example, the kernel \eqref{gamma-distr} defined by a gamma distribution
(see Section \ref{sect:gamma-dist} below), which
is approximated by
\begin{equation}
\label{k-approx}
k(t) \approx ch \sum_{n=M}^{N-1} \e^{\alpha nh} \e^{-(\e^{nh} + \kappa)t} .
\end{equation}
For large positive $n$ ($nh$ can be of size $50$ or more),
the exponent $\gamma_n = \e^{nh} + \kappa$ as well as
the coefficient $\e^{\alpha nh}$ are very large, and the computation has to be done
with care. For the kernel approximation above,
the integral term \eqref{distridelay} is approximated as
\begin{eqnarray} \nonumber
I (y)(t) & = & \int_0^t k(t-s)g\bigl( s, y(s)\bigr) \,\d s 
\approx  ch
\sum_{n=M}^{N-1} \e^{\alpha nh} z_n(t),\\[2mm] 
z_n(t) & \approx & \int_0^t  \e^{-\gamma_n (t-s)}g\bigl( s, y(s)\bigr) \,\d s  ,
\label{Ity}
\end{eqnarray}
where $z_n(t)$ is the solution of the differential equation (see Section~\ref{sect:dgln})
\[
\dot z_n(t) = -\gamma_n z_n(t) + g\bigl( t,y(t)\bigr) , \qquad z_n(0) = 0 .
\]
It is discretized numerically by a time integrator, typically a
Runge--Kutta method. For example, the $\theta$-method,
gives the approximation $z_n^k \approx z_n(t_k)$ via
\[
z_n^{k+1} = z_n^k - \gamma_n \Delta t \Bigl(\theta z_n^{k+1} + (1-\theta )z_n^k\Bigr)
+ \Delta t g_n^{k+1}, \quad g_n^{k+1} = g \Bigl( t_k+ \theta \Delta t , 
\theta y_n^{k+1} + (1-\theta )y_n^k \Bigr) ,
\]
which can also be written as
\[
z_n^{k+1} = R(-\gamma_n \Delta t) z_n^k + \Delta t S(-\gamma_n \Delta t) g_n^{k+1},\quad
R(\mu ) = \frac{1+(1-\theta ) \mu}{1-\theta \mu},\quad 
S(\mu ) = \frac{1}{1-\theta \mu} .
\]
Solving this recursion and using $z_n^0 = 0$ yields
\begin{equation}\label{znksum}
z_n^{k+1} = \Delta t  \sum_{j=0}^k R(-\gamma_n \Delta t)^{j} S(-\gamma_n \Delta t) g_n^{k+1-j} .
\end{equation}
Assuming $| g\bigl( t, y(t)\bigr) | \le G$ to be bounded, and
$|R(-\gamma_n \Delta t)| < 1$ we obtain (with $\mu = -\gamma_n \Delta t$)
\begin{equation}\label{znkest}
|z_n^{k+1}| \le \Delta t \big| S(\mu )\big| \frac{1 - | R(\mu )|^{k+1}}{1 - | R(\mu )|} G \le
\frac {\Delta t \big| S(\mu )\big| G }{1 - | R(\mu )|} .
\end{equation}
Any other Runge--Kutta method will lead to a similar recursion, where $R(\mu )$ is its
stability function, and $S(\mu )$ is a rational function bounded by $\bigo (\mu^{-1})$.
This shows that for a stable approximation of the integral term
\eqref{Ity}, $A_0$-stability (i.e., $|R(\mu )| <1$
for $\mu <0$) is a necessary condition, and $|R(\infty )|< 1$ is recommended.
These conditions imply that $|z_n^{k+1}| \le \bigo ( \Delta t |\mu|^{-1}|) =
\bigo ( \gamma_n^{-1})$.

Let us go back to the question of large coefficients in the approximation
\eqref{Ity}. In the exact solution,
a large coefficient $\e^{\alpha nh}$ is compensated by the small
factor $\e^{-(\e^{nh} + \kappa)t}$ in \eqref{k-approx}. In the numerical
approximation, it is compensated by $\gamma_n^{-1} \le \e^{-nh}$.

\subsection{Quasipolynomial case}

Consider now the situation, where some of the polynomials in \eqref{sum-exp-pol} have a
degree at least $1$. In this case $z_n(t)$ is a vector and satisfies
the differential equation
\begin{equation} \nonumber
\dot z_n(t) = J_n z_n(t) + e_n g\bigl( t , y(t)\bigr), \qquad z_n(0) = 0 \in \R^{m_n+1}.     
\end{equation}
A Runge-Kutta discretization yields
\begin{eqnarray} \nonumber
z_n^{k+1} & = & R(\Delta t J_n) z_n^k + \Delta t S(\Delta t J_n) e_n g_n^{k+1}\quad\text{and}\quad
\\
z_n^{k+1} & = & \Delta t  \sum_{j=0}^k R(\Delta t J_n)^{j} S(\Delta t J_n) e_n g_n^{k+1-j} 
\label{znkvecsum}
\end{eqnarray}
with $R(\Delta t J_n)$ and $S(\Delta t J_n)$ now rational matrix functions.
They are of the form (with $\gamma = \gamma_n$, $m=m_n$, and $\mu = - \gamma_n \Delta t$)
\[
\hskip -2mm R(\Delta t J_n) \!=\! \begin{pmatrix}
R(\mu )  & 0 & 0 &\cdots &  0 \\[-0.5mm]
{1 \choose 1} \Delta t R'(\mu) & R(\mu)  &0 & &\vdots \\[-0.5mm]
{2 \choose 2} \Delta t^2 R''(\mu) & {2 \choose 1} \Delta t R'(\mu) & R(\mu) & \ddots &  \vdots \\[-0.5mm]
\vdots & \vdots & \ddots & \ddots & 0 \\[0.5mm]
{m \choose m} \Delta t^m R^{(m)}(\mu) &{m \choose m-1} \Delta t^{m-1}R^{(m-1)}(\mu) & \cdots 
& {m \choose 1} \Delta t R'(\mu) & R(\mu)
\end{pmatrix} 
\]
This is a consequence of $ R(-\Delta t \gamma I + \Delta N) = 
\sum_{l= 0}^{m-1} \frac {\Delta t^l}{l!} R^{(l)}(\mu ) N^l$, where $N$ is the nilpotent matrix whose
only non-zero elements are in the first subdiagonal. For the matrix $S(\Delta t J_n)$
we get the same formulas with $R(\mu )$ replaced by $S(\mu )$.
We are interested to bound the norm of the vector $z_n^{k+1}$ in \eqref{znkvecsum}.
The first component $z_{n,0}^{k+1}$ is given by the relation \eqref{znksum} and therefore
satisfies the estimate \eqref{znkest}. To estimate the second component $z_{n,1}^{k+1}$
we have to compute the elements of the second row of $ R(\Delta t J_n)^{j}$.
The diagonal element is $R(\mu )^{j}$ and the element left to it is
$j\Delta t R' (\mu ) R(\mu )^{j-1}$.
From \eqref{znkvecsum} we thus get
\[
z_{n,1}^{k+1} = \Delta t \sum_{j=0}^k \Bigl( j \Delta t R' (\mu ) 
R(\mu )^{j-1} S(\mu ) + 
R(\mu )^{j} \Delta t S'(\mu ) \Bigr) g_n^{k+1-j} .
\]
Applying the triangle inequality, using the bound $| g\bigl( t, y(t)\bigr) | \le G$
and the estimates (for $0\le r <1$)
\begin{equation}\label{geomdiff}
\sum_{j=0}^k r^{j} = \frac{1-r^{k+1}}{1-r} \le \frac 1 {1-r}, \quad
\sum_{j=0}^k j\, r^{j-1}  \le \frac 1 {(1-r)^2} ,\quad
\sum_{j=0}^k {j \choose l}\, r^{j-l}  \le \frac 1 {(1-r)^{l+1}},
\end{equation}
(the second and last inequalities are obtained by differentiating the infinite geometric series
and by truncating higher order terms, which are positive)
we get
\begin{equation}\label{zn1kest}
|z_{n,1}^{k+1}| \le \Delta t^2 \biggl( \frac {\big| R' ( \mu ) 
S(\mu )\big|}
{(1 - |R (\mu )|)^2} + \frac {\big| S'(\mu )\big|}
{1-|R (\mu)|} \biggr) G .
\end{equation}
Applying the fact that $|S(\mu ) | = \bigo (\mu^{-1})$, $|S'(\mu ) | = \bigo (\mu^{-2})$,
and $|R'(\mu ) | = \bigo (\mu^{-2})$, we find that
$|z_{n,1}^k | = \bigo( \gamma_n^{-2})$, and the same conclusion as in Section~\ref{sect:stable}
can be drawn.

For general $m$ we use the multinomial theorem (the indices $i_l$ are non-negative)
\begin{equation}\label{multinomial}
\biggl(\, \sum_{l=0}^{m-1} \frac {{\Delta t}^l} {l!} \, R^{(l)}(\mu ) \,N^l \biggr)^j =
\sum_{i_1 + \ldots + i_m = j} {j \choose i_1, \ldots , i_m} \prod_{l=0}^{m-1}
\Bigl( \,\frac {\Delta t^l} {l!} \,  R^{(l)}(\mu ) \,N^l  \Bigr)^{i_{l+1}}
\end{equation}
and the inequalities \eqref{geomdiff} for values of $l$ up to $m$.
Since the nilpotent matrix satisfies $N^m = 0$, only the terms with
\begin{equation}\label{ijwsum}
i_2 + 2\,i_3 + \ldots + (m-1)\, i_m \le m-1
\end{equation}
give raise to non-vanishing expressions in \eqref{multinomial}. Therefore, only $i_1$ can be large.
It can take only the values $j,\,j-1, \ldots , j-m+1$, because otherwise
$i_2+ \ldots + i_m$ would have to be larger than $m-1$, contradicting the restriction
\eqref{ijwsum}.

We start with $i_1=j$. In this case all other $i_l$ are zero and there is only one
term in \eqref{multinomial}. The sum over $j$ from $0$ to $k$ is bounded by the
truncated geometric series (see \eqref{znkest}).

For $i_1=j-1$ only one index among $\{ i_2,\ldots ,i_m\}$ can be non-zero (equal to $1$,
if $m\ge 2$).
In each case the multinomial coefficient in \eqref{multinomial} equals $j$, and the
second sum in \eqref{geomdiff} provides the desired bound (see \eqref{zn1kest}).

For $i_1 = j-2$ we have two possibilities. Either one among the indices $\{ i_2,\ldots ,i_m\}$
equals $2$ (if $m\ge 3$), or two of them equal $1$ (provided that $m\ge 4$).
In each case the multinomial coefficient
in \eqref{multinomial} equals $j(j-1)/2$ or $j(j-1)$, and the estimate of
\eqref{geomdiff} with $l=2$ can be applied. Using $S''(\mu) = \bigo (\mu^{-3})$ and similar
estimates for the stability function $R(\mu)$, we obtain $|z_{n,2}^k | =\bigo ( \gamma_n^{-3})$.

This procedure can be continued until
$i_1 = j-m+1$. We obtain the bounds $|z_{n,l}^k | = \bigo ( \gamma_n^{-l-1})$ for the
$(l+1)$th component of the vector $z_n^k$.

\subsection{Determination of the accuracy parameters}
\label{sect:tol}

The variables $z_{i,j}(t)$ are auxiliary variables in the system \eqref{problem-aug} and the
question arises, how accurate they should be. For the original problem \eqref{problem}
the user has to specify the desired accuracy of the solution. This is typically done
with help of the parameters ${\it Atol}$ and ${\it Rtol}$, with the aim of having a
local error for $y(t)$ that is bounded by ${\it Atol} + |y| {\it Rtol}$  (for more details see
\cite[p.\ 124]{hairer96sod}). Of course, such a requirement can be applied component-wise.
For solution components with
a very small modulus that have a strong impact on the other solution components
it is in general advisable to use a
parameter  ${\it Atol}$ that is much smaller than ${\it Rtol}$. 

In the situation of the present work, it is not necessary that all variables $z_{i,j}(t)$
are computed very accurately. Since the solution $y(t)$ only depends on a linear
combination of them (with possibly very small coefficients), it is sufficient that this
linear combination is sufficiently accurate. This can be achieved as follows:
we augment the dimension of $y(t)$ by one and introduce the new variable $y_{d+1} (t)$ by
\begin{equation}\label{newvariable}
0 = \sum_{i=1}^n \sum_{j=0}^{m_i} c_{i,j} z_{i,j}(t) - y_{d+1} (t)  .
\end{equation}
We replace the double sum in \eqref{problem-aug} by $y_{d+1} (t)$ and we add the relation
\eqref{newvariable} to \eqref{problem-aug}. This gives a differential-algebraic system
where the number of algebraic variables is augmented by one. The advantage of this is that
we can require different accuracies for the $z_{i,j}$ and for their sum \eqref{newvariable}.
For a given tolerance {\it Tol} we propose (in general) to put
${\it Atol} = {\it Rtol} = {\it Tol}$ for all $d+1$ components of the augmented vector $y$,
and we put ${\it Atol} = {\it Rtol} = \omega {\it Tol}$ (with
$\omega \ge 1$) for all variables $z_{i,j}$. To increase efficiency without spoiling
accuracy we propose to take $\omega = 100$ or even larger
(see the second experiment in Section~\ref{sect:example1} and Table~\ref{tab:toldelta}).

\section{Approximating the kernel with a sum of exponentials}
\label{sect:sumofexp}

Important kernels, such as the gamma distribution and the Pareto distribution, contain
the function $t^{-\alpha}$ as factor. This section is devoted to approximate it by a sum
of exponential functions.

\subsection{Approach of Beylkin and Monz\' on}
\label{sect:beylkin}

Since the Laplace transform of the function $z^{\alpha -1}$ is $\Gamma (\alpha )\, t^{-\alpha }$
for $\alpha >0$,
we have the integral re\-pre\-sen\-ta\-ti\-on
\begin{equation}\label{integral-repr}
t^{-\alpha} = \frac1 {\Gamma (\alpha )} \int_0^\infty \e^{-tz}z^{\alpha -1} \,\d z =
 \frac1 {\Gamma (\alpha )} \int_{-\infty}^\infty \e^{-t \e^s}\e^{\alpha  s} \,\d s ,\qquad \alpha > 0 .
\end{equation}
To express this function as a sum of exponentials, \cite{beylkin10abe} proposes to
approximate the integral to the right by the trapezoidal rule. With a step size $h>0$ this yields
\begin{equation}\label{Tth}
T (t,h) = \frac h {\Gamma (\alpha )} \sum_{n=-\infty }^\infty \e^{\alpha nh} \e^{-\e^{nh} t} .
\end{equation}

\noindent
{$\bullet$\quad\it Error of the trapezoidal rule.}
It follows from \cite[Theorem 5.1]{trefethen14tec} that the error due to this approximation
can be bounded by
\[
\big| t^{-\alpha } - T(t,h) \big| \le 
\frac {2C}{\e^{2\pi a / h} -1}, \qquad C = \sup_{b\in (-a,a)} 
\frac 1 {\Gamma (\alpha )}\int_{-\infty}^\infty 
\e^{-t \e^s \cos b } \e^{\alpha s} \,\d s   = (t \cos a )^{-\alpha} 
\]
for any $h>0$ and any $0<a< \pi /2$. This gives the estimate
\begin{equation}\label{est-T}
\big| t^{-\alpha } - T(t,h) \big| \le c_\alpha t^{-\alpha}, \qquad
c_\alpha = \frac{2 (\cos a)^{-\alpha}}{\e^{2\pi a /h} -1} .
\end{equation}
To achieve $c_\alpha \le \eps$, the step size $h$ has to satisfy the inequality in
\begin{equation}\label{stepsize}
h\le \frac{2\pi a}{\ln \bigl( 1 + \frac 2\eps (\cos a )^{-\alpha} \bigr) } , \qquad 
a = \frac \pi 2 \Big( 1 - \frac \alpha {(\alpha +1)\, \ln \eps^{-1}} \Big) .
\end{equation}
The value of $a$ is chosen so that, for
a given $\eps >0$, the estimate for $h$ in \eqref{stepsize} turns out to be close to maximal.

\medskip
\noindent
{$\bullet$\quad\it Error from truncating the series \eqref{Tth} at $ - \infty$.}
Since the function $s\mapsto \e^{-t \e^s}\e^{\alpha  s}$ is monotonically increasing
for $s\le \ln (\alpha /t)$ we have, for $M h\le \ln (\alpha /t)$,
\begin{equation}\label{TthM}
\begin{array}{rcl}
 E_M (t,h) & = & \displaystyle
 \frac h {\Gamma (\alpha )} \sum_{n=-\infty }^{M-1} \! \e^{\alpha nh} \e^{-\e^{nh} t} 
\le \displaystyle  \frac1 {\Gamma (\alpha )} \int_{-\infty}^{Mh} \!\! \e^{-t \e^s}\e^{\alpha  s} \,\d s \\[4mm]
&=&
\displaystyle \frac{t^{-\alpha}}{\Gamma (\alpha )} \int_0^{t \e^{Mh}} \!\!\!\!\e^{-\sigma } 
\sigma^{\alpha -1}\, \d \sigma
= t^{-\alpha } \biggl(1 - \frac{\Gamma (\alpha ,t \e^{Mh})}{\Gamma (\alpha )}\biggr) ,
\end{array}
\end{equation}
where $\,\Gamma (\alpha , x) = \int_x^\infty \e^{-\sigma } 
\sigma^{\alpha -1}\, \d \sigma \,$ is the incomplete Gamma function.

\medskip
\noindent
{$\bullet$\quad\it Error from truncating the series \eqref{Tth} at $ + \infty$.}
The function $s\mapsto \e^{-t \e^s}\e^{\alpha  s}$ is monotonically decreasing
for $s\ge \ln (\alpha /t)$.
The same computation as before therefore shows that for $Nh \ge \ln (\alpha /t)$
\begin{equation}\label{TthN}
\begin{array}{rcl}
E^N (t,h) & = & \displaystyle
\frac h {\Gamma (\alpha )} \sum_{n=N }^{\infty} \e^{\alpha nh} \e^{-\e^{nh} t} 
\le  \frac1 {\Gamma (\alpha )} \int_{Nh}^{\infty} \e^{-t \e^s}\e^{\alpha  s} \,\d s 
\\[4mm]
& = & \displaystyle
\frac{t^{-\alpha}}{\Gamma (\alpha )} \int_{t \e^{Nh}}^\infty \e^{-\sigma } \sigma^{\alpha -1}\, \d \sigma
= t^{-\alpha}\,\frac{\Gamma (\alpha ,t \e^{Nh})}{\Gamma (\alpha )} .
\end{array}
\end{equation}

\begin{figure}[h!]
\centering
 \begin{picture}(0,0)
  \epsfig{file=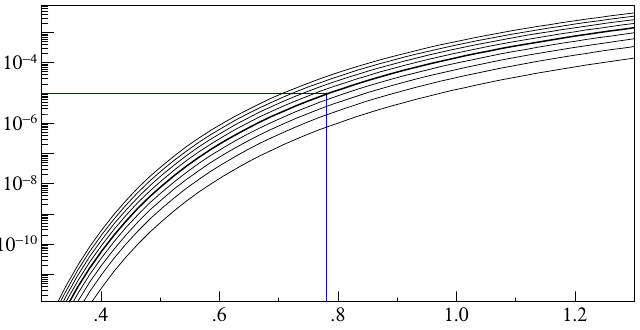}
 \end{picture}%
\begin{picture}(308.6,157.8)( -21.3, 442.4)
  \GGGput[  -15,  140](  0.00,597.43){error of trapezoidal rule}
  \GGGput[    1,    0](150.68,536.24){$\alpha = 0.1$}
  \GGGput[  100,    0](142.25,569.01){$\alpha = 0.9$}
  \GGGput[  130,  -60](284.50,455.18){step size $h$}
 \end{picture}
\vspace{-3mm}
\caption{Relative error
$\big| t^{-\alpha } - T(t,h) \big|\big/ t^{-\alpha}$
as a function of the step size $h$. The different curves correspond to
$\alpha = 0.1, 0.2, \ldots , 0.9$ (bold curve for $\alpha = 0.5$).
\label{fig:err-trap}}
\end{figure}

\medskip
\noindent
{$\bullet$\quad\it Choice of the step size $h$.}
For a fixed accuracy requirement $\eps$,
we choose $h$ according to \eqref{stepsize}.
To study the quality of this approximation, we plot
in Figure~\ref{fig:err-trap} for several values of $\alpha$
the relative error of the trapezoidal rule as a function of the step size $h$.
We see that for $\alpha =1/2$ the error of the trapezoidal rule is bounded
by $\eps = 10^{-5}$ whenever $h\le 0.78$. The estimate \eqref{stepsize} gives
$h\le 0.70$, which is a reasonably good approximation.

\medskip
\noindent
{$\bullet$\quad\it Choice of the truncation indices $M$ and $N$.}
For a treatment of the two truncation errors, we restrict $t$ to the interval
$0<\delta \le t \le T < \infty$, and we use the estimates \eqref{TthM} and \eqref{TthN}.
The incomplete Gamma function is monotonically decreasing and satisfies
$\Gamma (\alpha ,0 ) = \Gamma (\alpha )$, $\Gamma (\alpha ,\infty ) = 0$.
We let $x_*$ be such that $\Gamma (\alpha , x_*) \ge \Gamma (\alpha ) ( 1- \eps )$, and
$x^*$ be such that $\Gamma (\alpha , x^*) \le \Gamma (\alpha )\, \eps $.
Both truncation errors are bounded by $\eps\, t^{-\alpha} $ for $t\in [\delta , T ]$ if
$M$ and $N$ are chosen according to $T\e^{Mh} \le x_*$ and
$\delta \e^{Nh}  \ge x^*$.
Since $\Gamma (\alpha ,x) \ge \Gamma (\alpha ) - x^\alpha / \alpha$,
we can approximate $x_*$ by the relation $x_* ^\alpha =  \Gamma (\alpha +1) \, \eps $.
From $ \Gamma (\alpha ,x ) \le x^{\alpha -1}\e^{-x}$ a suitable $x^*$
is given by $ (x^*)^{\alpha -1}\e^{-x^*} \le \Gamma (\alpha )\, \eps$, which is approximately
$x^* = - \ln \bigl( \Gamma (\alpha )\, \eps \bigr) $.

\begin{figure}[h]
\centering
 \begin{picture}(0,0)
  \epsfig{file=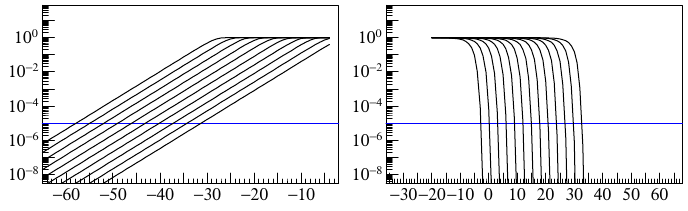}
 \end{picture}%
\begin{picture}(331.4,100.8)( -10.0, 499.4)
  \GGGput[  -10,  140]( 11.32,597.43){\small truncation error at $-\infty$}
  \GGGput[  190, -100](153.57,512.15){\small$M$}
  \GGGput[    0,    0](115.15,548.77){\small$T=10$}
  \GGGput[    0,    0]( 21.32,563.34){\small$T=10^{10}$}
  \GGGput[  -10,  140](176.34,597.43){\small truncation error at $+\infty$}
  \GGGput[  190, -100](318.59,512.15){\small$N$}
  \GGGput[    0,    0](270.29,571.90){\small$\delta =10^{-10}$}
  \GGGput[    0,    0](183.45,524.92){\small$\delta=10^2$}
 \end{picture}
\vspace{-3mm}
\caption{Relative error of the the truncation errors as a function of $M$ and $N$, respectively.
\label{fig:err-trunc}}
\end{figure}

To get an impression of the size of the
truncation indices $M$ and $N$, we fix the values $\alpha =0.5$ and $h=0.78$, which correspond to
$\eps = 10^{-5}$
(see Figure~\ref{fig:err-trap}).
In the left picture of Figure~\ref{fig:err-trunc} we plot 
$\max \big\{ E_M(t,h)/t^{-\alpha} \big|  t\in[10^{-10},T] \big\}$ (the relative truncation error) 
as a function of $M$,
for $T=10, 10^2, \ldots , 10^{10}$. 
In the right picture of Figure~\ref{fig:err-trunc} we plot 
$\max \big\{ E^N(t,h)/t^{-\alpha} \big|  t\in[\delta , 500 ] \big\}$ 
as a function of $N$,
for $\delta =10^2, 10, 1, 10^{-1}, \ldots , 10^{-10}$. Replacing the right end
of the interval $[\delta , 500]$ by a larger value does not change the picture.
We have included the curves for large values of $\delta$, because they are of interest
for a treatment of the Pareto distribution (see Section~\ref{sect:pareto}).

As a conclusion we see that for an accuracy requirement of $\eps = 10^{-5}$ and for
an interval $[\delta , T]$ with $\delta = \eps = 10^{-5}$ and $T=100$, we can choose
the step size $h=0.78$, and the truncation indices $M=-35$ and $N=18$ or
(using $h$ from \eqref{stepsize})
$h=0.70$, $M=-40$ and $N=20$. 

\medskip
\noindent
{$\bullet$\quad\it Total approximation error.}
The error of approximating the function $t^{-\alpha}$ by a sum of exponentials is composed by
the error of the trapezoidal rule and by those due to truncation. We have,
for $\delta \le t \le T$,
\begin{equation}\label{total-err}
\big| t^{-\alpha} - T_M^N (t,h) \big| \le 3\, \eps \, t^{-\alpha}, \qquad 
T_M^N (t,h) = \frac h {\Gamma (\alpha )} \sum_{n=M }^{N-1} \e^{\alpha nh} \e^{-\e^{nh} t} ,
\end{equation}
provided that $h$, $M$, and $N$ are chosen as discussed above.

\subsection{Approach of Braess and Hackbusch}
\label{sect:_braess}

An efficient computation of high-dimensional integrals is another motivation for
studying the approximation of $t^{-\alpha}$ ($\alpha = 1/2$ or $1$) by exponential sums
(see \cite{braess09ote}). The differences to the approach of Beylkin and Monz\'on
are that the coefficients $c_i$ and $\gamma_i$ in the approximation
$u_n(t) = \sum_{i=1}^n c_i \exp{(-\gamma_i t)}$
are computed for minimal $n$ numerically with a Remez-like algorithm, and the absolute error
$\big| t^{-\alpha} - u_n(t) \big|$ (not the relative error as in \eqref{total-err}) is kept below a tolerance $\eps$.

\section{Important examples of distributed delays}
\label{sect:distributed}

Because of their application in pharmacodynamics we are mainly interested in the
situation where the kernel of the integral delay term \eqref{distridelay} is a distribution,
i.e., a non-negative
function with $\int_0^{\infty} k(s)\, \d s =1$.
The following distributions are taken from the list given in
\cite{wesolowski20cot}.

\subsection{Gamma distribution}
\label{sect:gamma-dist}

For positive parameters $\kappa > 0$ and $0 < \alpha < 1$, the gamma distribution is given by
\begin{equation}\label{gamma-distr}
k(t) = \frac {\kappa^{1-\alpha}}{\Gamma (1-\alpha)}\, t^{-\alpha} \e^{-\kappa t}, \qquad
t > 0 .
\end{equation}
We choose $\delta > 0$ such that $\int_0^\delta k(t)\, \d t \le \eps$, and $T>0$ such that
$\int_T^\infty k(t)\, \d t \le \eps$. Roughly speaking, the small parameter $\delta$
and the large
parameter $T$ have to be such that 
\begin{equation} \label{eq:TGamma}
\frac{(\kappa \,\delta  )^{1-\alpha} }{\Gamma (2-\alpha )} \le \eps\qquad \text{and} \qquad
\frac{(\kappa \,T)^{-\alpha} \e^{-\kappa T} }{\Gamma (1-\alpha )} \le \eps.
\end{equation}
We use \eqref{total-err} for approximating the gamma distribution by a sum of exponential
functions.
For a given accuracy requirement $\eps$,
Algorithm~\ref{alg_paramG} summarizes the computation of the parameters
$\delta , T$ as well as the parameters $ h, M, N$ needed in \eqref{total-err}.

\medskip
\begin{algorithm}[H] \label{alg_paramG}
\DontPrintSemicolon
\KwData{$\alpha$, $\kappa$, $\eps$, $\delta_{\min}$, $t_f$}
\KwResult{$h, T, \delta, M, N$} 
\Begin{
\nl Set $a = \displaystyle \frac \pi 2 \Big( 1 - \frac \alpha {(\alpha +1)\, \ln \eps^{-1}} \Big)$\;
\nl Set $h = \displaystyle \frac{2\pi a}{\ln \bigl( 1 + \frac 2\eps (\cos a )^{-\alpha} \bigr) }$\;
\nl Compute $T$ such that $ (\kappa \,T)^{-\alpha} \e^{-\kappa T} / \,\Gamma (1-\alpha ) = \eps$\;
\nl Set $T = \min\{ t_f, T \}$\;
\nl Set $x_* =  \left(\Gamma (\alpha +1) \, \eps \right)^{1/\alpha}$\;
\nl Set $x^* = - \ln \bigl( \Gamma (\alpha )\, \eps \bigr)$\;
\nl Compute $M$ such that $T\e^{M h} = x_*$\;
Compute $\delta$ such that $(\kappa \,\delta  )^{1-\alpha} / \,\Gamma (2-\alpha ) = \eps$\;
\nl Set $\delta = \max\{ \delta, \delta_{\min} \}$\;
\nl Compute $N$ such that $\delta \e^{N h}  = x^*$\;
\Return
}
\caption{\rule{0mm}{3mm}Choice of the parameters for the Gamma distribution.}
\end{algorithm} 

\medskip
\noindent
{$\bullet$\quad\it Negative values of $\alpha$.} 
For negative values of $\alpha$ the function $k(t)$ vanishes at the origin, but the sum of
exponentials \eqref{Tth} does not. Let $\alpha$ satisfy $-1<\alpha < 0$. The idea is to split
$t^{-\alpha} = t \cdot t^{-\alpha -1}$, so that $\alpha +1 \in (0,1)$, and the formula
\eqref{Tth} can be applied to the factor $t^{-\alpha -1 }$.
This results in a linear combination of functions $t \,\e^{-\gamma t}$, which can be
treated as explained in Section~\ref{sect:dgln}.

For values of $\alpha$ satisfying $-2< \alpha < -1$, we split
$t^{-\alpha } = t^2\cdot t^{-\alpha -2}$.
The above procedure then leads to kernels of the form \eqref{sum-exp-pol} with $m_i=2$.
An extension
to other non-integer negative values of $\alpha$ is straight-forward.

\subsection{Pareto distribution}
\label{sect:pareto}

The so-caled 
type I Pareto distribution \cite{wesolowski20cot} is given for $\alpha >0$ and $\beta >0$ by
\begin{equation}\label{pareto-distr}
k(t) = \Bigg\{\begin{array}{cl} 0 & 0\le t < \beta \\[2mm]
\alpha \,\beta^\alpha\, t^{-\alpha -1}~~~  & t \ge \beta .\end{array}
\end{equation}
Note that there is no singularity at the origin. For the
approximation \eqref{total-err} of $t^{-\alpha -1}$, we can choose $\delta = \beta$ independent
of the accuracy $\eps$, so that the condition $\delta \e^{Nh} \ge x^*$ is less restrictive for the truncation index $N$. The condition $\int_T^\infty k(t)\, \d t \le \eps$ for $T$ reduces to
$T^{\alpha} \ge \beta^\alpha  / \eps$. The truncation index $M$ is determined by
$T \e^{Mh} \le x_*$, where $T$ is the minimum of $\beta /\eps^{1/\alpha}$ and
of the length of the integration interval $T_{end}$
(see Algorithm~\ref{alg_paramP}).

\begin{algorithm}[ht] \label{alg_paramP}
\DontPrintSemicolon
\KwData{$\alpha$, $\beta$, $\eps$, $t_f$}
\KwResult{$h, T, M, N$} 
\Begin{
\nl Set $T = \beta\,\eps^{-1/\alpha}$\;
\nl Set $T = \min\{ t_f, T \}$\;
\nl Set $a = \displaystyle \frac \pi 2 \Big( 1 - \frac {(\alpha+1)} {(\alpha + 2)\, \ln \eps^{-1}} \Big)$\;
\nl Set $h = \displaystyle \frac{2\pi a}{\ln \bigl( 1 +  2\eps^{-1} (\cos a )^{-(\alpha+1)} \bigr) }$\;
\nl Set $x_* =  \Gamma (\alpha +2) \, \eps $\;
\nl Set $x^* = - \ln \bigl( \Gamma (\alpha+1 )\, \eps \bigr)$\;
\nl Compute $M$ such that $T\e^{M h} = x_*$\;
Compute $N$ such that $\beta \e^{N h}  = x^*$\;
\Return
}
\caption{\rule{0mm}{3mm}Choice of the parameters for the Pareto distribution.}
\end{algorithm} 

In this case the integral in \eqref{distridelay} is from $0$ to $t-\beta$. 
When replacing the kernel (on the interval $[\beta , T]$) by a sum of exponential functions, 
we are thus lead to
\begin{equation}\label{eqyparet}
M \dot y (t) = f \bigl( t, y(t), y(t-\tau ), I_P ( t - \beta ) \bigr)
\end{equation}
where the function $I_P (t)$ is given by (with $\gamma_n = \e^{nh}$)
\begin{equation} \label{def-Ip}
I_P (t) = \left\{ \begin{array}{cr} 0 & \mbox{if\quad$t \le 0$} \\[2mm]
\displaystyle \alpha \beta^\alpha \frac{h}{\Gamma (\alpha )}
\sum_{n=M}^{N-1} \e^{(\alpha +1)nh} \e^{-\gamma_n \beta} z_n(t ) ~~& 
\mbox{if\quad$t \ge 0$}
\end{array}
\right.
\end{equation}
and $\displaystyle
~ z_n(t) = \int_0^t \e^{-\gamma_n (t-s)} g\bigl( s, y(s)\bigr) \, \d s ~  $
gives raise
to the differential equation
\begin{equation}\label{znparet}
\dot z_n (t) = - \gamma_n z_n(t) +  g\bigl( t , y(t ) \bigr),\qquad
z_n(0) = 0 .
\end{equation}
As we have done in Section~\ref{sect:tol} we consider $I_p(t)$ as a new variable
of the system, and we add the algebraic relation \eqref{def-Ip} to the
augmented system \eqref{eqyparet}-\eqref{znparet}. In this way we can
monitor the choice of the accuracy parameters ${\it Atol}$ and
${\it Rtol}$ as explained in Section~\ref{sect:tol}. The introduction of the
variable $I_P (t)$ in the system has the further advantage that during the
numerical integration only the back values of the scalar function
$I_P(t)$ have to be stored.
Note that the augmented system has additional
breaking points at $i\tau + j \beta$ with non-negative integers $i,j$.

\medskip
\noindent
{$\bullet$\quad\it Remark.} For small positive $\alpha$ (e.g., $\alpha = 0.15$ and smaller, see
\cite[Table 1]{wesolowski20cot}) the condition $T^{\alpha} \ge \beta^\alpha  / \eps$ can lead
to a very large $T$, which then results in a large negative truncation index $M$. To increase
efficiency, we propose to introduce the new variable
\[
y_{d+1}(t) = \alpha \beta^\alpha
\int_0^{t-\beta} (t-s)^{-\alpha -1} g\bigl( s, y(s)\bigr) \, \d s ,
\]
which is the argument $I(y)(t)$ in \eqref{problem},
and to differentiate it with respect to $t$. This yields the delay differential equation
\[
\dot y_{d+1}(t) = \alpha \beta^{-1} g\bigl( t-\beta , y(t-\beta ) \bigr) -
\alpha (\alpha + 1)\beta^\alpha
\int_0^{t-\beta} (t-s)^{-\alpha -2} g\bigl( s, y(s)\bigr) \, \d s ,
\]
which can be added to the system for $y(t)$. The integral in the right-hand side is again of
Pareto type, but with $\alpha$ replaced by $\alpha + 1$. Replacing the new kernel
by a sum of exponentials  gives a more efficient algorithm. Of course, this procedure
can be repeated to increase the parameter $\alpha$ even more.

\section{Numerical experiments}
\label{sect:numerical}

This section provides numerical evidence of the proposed algorithm. The first two
examples are scalar test equations (one ordinary differential equation with a
Gamma-distributed delay term, the other a delay differential equation with a
Pareto-distributed delay term). The third example is taken from applications in
chemotherapy-induced myelo-suppression.
For the time integration we make use of the code {\sc Radau5}
(see\cite{hairer96sod}), if the augmented
system is an ordinary differential equation, and of the
code {\sc Radar5} (see \cite{guglielmi01irm,guglielmi08cbp}), if the augmented system
is a delay differential equation.

\subsection{Example 1: ordinary differential equation with a Gam\-ma-di\-stri\-buted delay term}
\label{sect:example1}

As a first test example we consider a linear differential equation
\begin{equation}\label{example2}
\dot y (t) =   \bigl(1 - y(t) \bigr)\,\text{erf}\left(\frac{\sqrt{t}}{2}\right) -
   \frac{\e^{-t/4}\sqrt{t}}{\sqrt{\pi }}  + I\bigl( t,y(t) \bigr) + \frac12,
\qquad y(0)=0,
\end{equation}
where $y(t)$ is a real-valued function, $\text{erf}$ denotes the error function,
and the distributed delay term is 
\begin{equation}\label{distridelay2}
I (y)(t) = \int_0^t k(t-s) y(s) \,\d s , \qquad k(t) = \frac{\e^{-t/4}}
{2 \sqrt{\pi t}}
\end{equation}
Here, $k(t)$ is the Gamma-distribution \eqref{gamma-distr} with parameters
$\kappa = 1/4$ and $\alpha =1/2$.
Note that the kernel is weakly singular, but that there is no singularity in
the augmented system \eqref{problem-aug}.
The inhomogeneity of the equation is chosen such that
\[
y(t)=t/2
\]
is the exact solution. This follows from the fact that with $y(t)=t/2$ we have
\[
I (y)(t) = \frac{t-2}{2}
   \,\text{erf}\left(\frac{\sqrt{t}}{2}\right)+
   \frac{e^{-t/4}\sqrt{t}}{\sqrt{\pi }}.
\]

\noindent{\it
First experiment (connection between ${\it Tol}$ and $\eps$).}
We fix the tolerance ${\it Tol} = 10^{-8}$
for the time integrator, and we vary
the value of $\eps$, which determines the accuracy of the approximation by the sum
of exponentials.
We consider the integration interval $[0,t_f]$ with $t_f=50$, so that the exact solution
at the final point is $y(t_f)= 25$.

According to Algorithm~\ref{alg_paramG} of Section~\ref{sect:gamma-dist} we compute,
for given $\eps >0$ and for $\delta_{\min} =0$, the parameters $h$, $T$, $M$, and $N$.
For $\eps = {\it Tol} =10^{-8}$, we compute $T=\min\{50,65.79\} = 50$ and
$\delta = \pi \cdot 10^{-16}$.
The suggested parameter values are $h=0.4638$, $M=-89$, $N=84$ which are
indicated in bold in Table \ref{tab:21}. There, the values of $h,T,M,N$ are also given
for further values of $\eps$. The value of $\delta$ is $\delta = \pi \eps^2$.

\begin{table}[ht]
\begin{center}
{\small
\begin{tabular}{|c||c|c|c|c|c|} \hline
 $\eps$ & $h$ & $T$ & $M$ & $N$ & ${\rm err}$ \\ 
\hline
\hline
\rule[0mm]{0cm}{3.5mm}%
$10^{-4}$  &  $ 0.84$  &  $30.49$  &  $ -27$  &  $  24$  &  $2.45 \cdot 10^{-4}$ \\
 $10^{-5}$  &  $ 0.70$  &  $39.20$  &  $ -39$  &  $  35$  &  $2.75 \cdot 10^{-5}$ \\
 $10^{-6}$  &  $ 0.60$  &  $48.00$  &  $ -54$  &  $  49$  &  $2.35 \cdot 10^{-6}$ \\
 $10^{-7}$  &  $ 0.52$  &  $50.00$  &  $ -70$  &  $  65$  &  $2.40 \cdot 10^{-7}$ \\
 $\mathbf{10^{-8}}$  &  $\mathbf{0.46}$  &  $\mathbf{50.00}$  &  $\mathbf{-89}$  &  
 $\mathbf{84}$  &  $\mathbf{1.71 \cdot 10^{-8}}$ \\
 $10^{-9}$  &  $ 0.42$  &  $50.00$  &  $-110$  &  $ 104$  &  $4.72 \cdot 10^{-10}$ \\
 $10^{-10}$  &  $ 0.38$  &  $50.00$  &  $-133$  &  $ 127$  &  $2.14 \cdot 10^{-9}$ \\
 $10^{-11}$  &  $ 0.35$  &  $50.00$  &  $-158$  &  $ 152$  &  $2.08 \cdot 10^{-9}$ \\
  \hline
\end{tabular}
}
\vskip 1.5mm
\caption{Error behavior obtained by applying {\sc Radau5} to the test problem \eqref{example2} with final point $t_f=50$ and ${\it Tol}= 10^{-8}$ for different computed values of
$h$, $M$, and $N$ (corresponding to $\eps$).}
\label{tab:21} 
\end{center}
\end{table} 

For the numerical integration we apply the code {\sc Radau5}
(the code {\sc Radar5} would give similar results)
with accuracy requirement ${\it Atol} = {\it Rtol} = {\it Tol} =10^{-8}$
to the augmented system \eqref{problem-aug} with parameters $h,M,N$, determined by
different values of $\eps$. The initial step size is $h=\eps$.
The last column of
Table \ref{tab:21} reports the relative error
$\,{\rm err} := |y(t_f) - \bar{y}|/|y(t_f)| \,$
at the final point.
Here, $\bar{y}$ denotes the numerical approximation computed at time $t_f$.
We can observe that for $\eps\le {\it Tol}$ the error is essentially proportional to
$\eps$, which corresponds to the error of the kernel approximation. For
$\eps \ge {\it Tol}$ the error remains close to ${\it Tol}$, which shows the
error of the time integration.

\begin{table}[h] 
\begin{center}
{\small
\begin{tabular}{|c||c|c|c||c|c|c||c|c|c|} 
\hline
&\multicolumn{3}{|c||}{\rule{0mm}{3mm}$\omega =1$}&\multicolumn{3}{|c||}{$\omega =10$}
&\multicolumn{3}{|c|}{$\omega =100$} \\
\hline
 $\eps$ & ${\rm err}$ & $\# {\rm fe}$ & ${\rm cpu}$ & 
${\rm err}$ & $\# {\rm fe}$ & ${\rm cpu}$ & 
${\rm err}$ & $\# {\rm fe}$ & ${\rm cpu}$ \\
\hline
\hline
\rule[0mm]{0cm}{3.5mm}%
$10^{-4}$  &  $  2.5\,e\text{-}4$  &  $ 81$  &  $ 2.5\,e\text{-}4$  & 
$2.5\,e\text{-}4$   &  $   66$  &  $  2.2\,e\text{-}4$  & 
$2.5\,e\text{-}4$  &  $    66$  &  $  2.2\,e\text{-}4$   \\
$10^{-6}$  &  $  2.4\,e\text{-}6$  &  $ 162$  &  $ 7.9\,e\text{-}4$  & 
$2.3\,e\text{-}6$   &  $   132$  &  $  6.6\,e\text{-}4$  & 
$2.3\,e\text{-}6$  &  $    117$  &  $  6.0\,e\text{-}4$   \\
$10^{-8}$  &  $  1.8\,e\text{-}8$  &  $ 365$  &  $ 2.7\,e\text{-}3$  & 
$1.6\,e\text{-}8$   &  $   279$  &  $  2.1\,e\text{-}3$  & 
$1.5\,e\text{-}8$  &  $    243$  &  $  1.8\,e\text{-}3$   \\
$10^{-10}\!\!$  &  $  5.8\,e\text{-}11\!\!$  &  $ 773$  &  $ 7.6\,e\text{-}3$  & 
$1.1\,e\text{-}11\!\!$   &  $   587$  &  $  6.2\,e\text{-}3$  & 
$1.2\,e\text{-}10\!\!$  &  $    482$  &  $  5.1\,e\text{-}3$   \\
\hline
\end{tabular}
}
\vskip 1.5mm
\caption{Error behavior of {\sc Radau5} applied to the test problem
\eqref{example2} with $t_f=50$, initial step size $h=0.1$, and
${\it Atol} = {\it Rtol} = \eps$ for $y$ and $y_{d+1}$
(see \eqref{newvariable}),
and ${\it Atol} = {\it Rtol} = \omega \eps$ for the auxiliary variables $z_{i,j}$.
Here, ${\rm err}$, $\# {\rm fe} $, $ {\rm cpu}$ indicate the relative error of $y(t)$,
the number of function evaluations and the cpu time.}
\label{tab:toldelta} 
\end{center}
\end{table} 

\medskip
\noindent{\it
Second experiment (choice of the accuracy parameters).}
As proposed in Section~\ref{sect:tol} we add a new variable $y_{1+1}(t)$ to the
equation \eqref{example2}, and we apply the code {\sc Radau5} to the augmented system
with accuracy parameters ${\it Atol} = {\it Rtol} = \eps$ for the $y$-components
and ${\it Atol} = {\it Rtol} = \omega \eps$ for the $z$-components. We first put
$\omega = 1$, and then study the behaviour for increasing $\omega$.
The relative error of $y(t)$, the number of function evaluations, and the cpu
time are given in Table~\ref{tab:toldelta}. We observe that increasing $\omega$ does not affect
the precision of the numerical approximation for $y(t)$, but significantly reduces the
number of function evaluations and the cpu time. This is due to the fact that
the numerical integrator can take larger step sizes.
For every $\eps$ there seems to exist a critical value of
$\omega$, such that the numerical result does not change any more when $\omega$ is
further increased. 
Numerical computations show that this happens for $\eps = 10^{-4}$ when $\omega = 10$,
for ${\eps}=10^{-6}$ when $\omega = 100$, for $\eps = 10^{-8}$ when 
$\omega = 1000$, etc.
We thus propose to use for this example $\omega = 0.1 \eps^{-1/2}$, i.e.,
${\it Rtol} = \sqrt{\eps}/10$.

\subsection{Example 2: delay differential equation with a Pareto-distributed delay term}

We next consider a test example with Pareto-distributed delay term.
Since such a kernel leads in any case to a delay equation, we add a
discrete delay term and consider
\begin{equation}\label{example1}
\dot y (t) =  - 5\,I (y)(t)  -
 \frac{y\bigl( t - \tau \bigr) - 2 }{y(t) + 1},
\qquad  y(t) = t ~~ \text{for} ~~ t\le 0 ,
\end{equation}
where $y(t)$ is a scalar real function, the discrete delay term is
$ \tau = \pi/4$, and
the distributed delay term is 
\begin{equation}\label{distridelay1}
I (y)(t) = \frac{\sqrt \beta}{2}\int_0^{t-\beta} {(t-s)^{-3/2}}\, y(s) \,\d s ,
\end{equation}
where the kernel is a Pareto distribution \eqref{pareto-distr} with $\alpha = 1/2$.
Since $t-s \ge \beta > 0$, there is no singularity in the integrand.
Note that the system \eqref{eqyparet} depends on $y(t-\tau)$ and also $I_P(t-\beta )$, which gives rise to
breaking points at integral multiples of $\tau$ and $\beta$ and also their integer linear combinations.

\begin{table}[th] 
\begin{center}
{\small
\begin{tabular}{|c||c||c|c|c|c|} \hline
$ \eps$ & $h$ & $M$ & $N$ & ${\rm err}$\\ 
\hline
\hline
\rule[0mm]{0cm}{3.5mm}%
 $10^{-1}$  &  $  1.662$  &  $  -3$  &  $   1$  &  $7.69 \ \cdot 10^{-2}$   \\
 $10^{-2}$  &  $  1.116$  &  $  -6$  &  $   2$  &  $8.97 \ \cdot 10^{-4}$   \\
 $10^{-3}$  &  $  0.851$  &  $ -11$  &  $   3$  &  $2.31 \ \cdot 10^{-4}$   \\
 $10^{-4}$  &  $  0.692$  &  $ -17$  &  $   4$  &  $2.81 \ \cdot 10^{-5}$   \\
 $10^{-5}$  &  $  0.586$  &  $ -24$  &  $   5$  &  $1.37 \ \cdot 10^{-6}$   \\
 $10^{-6}$  &  $  0.509$  &  $ -32$  &  $   6$  &  $3.46 \ \cdot 10^{-7}$   \\
 $10^{-7}$  &  $  0.451$  &  $ -41$  &  $   7$  &  $1.90 \ \cdot 10^{-7}$   \\
 $10^{-8}$  &  $  \mathbf{0.405}$  &  $ \mathbf{-51}$  &  $\mathbf{8}$  &  $\mathbf{9.83 \ \cdot 10^{-8}}$   \\
 $10^{-9}$  &  $  0.368$  &  $ -62$  &  $   9$  &  $5.95 \ \cdot 10^{-8}$   \\
 $10^{-10}$  &  $  0.337$  &  $ -75$  &  $  10$  &  $1.75 \ \cdot 10^{-7}$   \\
 $10^{-11}$  &  $  0.311$  &  $ -88$  &  $  11$  &  $2.40 \ \cdot 10^{-7}$   \\
  \hline
\end{tabular}
}
\vskip 1.5mm
\caption{Error behavior obtained by applying {\sc Radar5} to the test problem \eqref{example1}. for different values of
$h$, $M$ and $N$, computed according Algorithm \ref{alg_paramP},
with ${\it Tol}= 10^{-8}$}
\label{tab:1} 
\end{center}
\end{table}

In our experiments we set $\beta=1$ and we
consider the integration interval $[0,t_f]$ with $t_f=10$.
With this choice we compute a reference solution to high precision,
\[
y(t_f) \approx 0.570525788119.
\]

As explained in Section~\ref{sect:pareto}, an approximation of the kernel by a
sum of exponentials gives rise to a delay differential equation. For an approximation
error $\eps$ Algorithm~\ref{alg_paramP} yields the parameters $h$, $M$, and $N$
that are needed for the description of the augmented system.

For $\eps = 10^{-8}$, the values are $h=0.405$, $M=-51$ and $N=8$, which are
indicated in bold in the right side of Table~\ref{tab:1}. 
In contrast to the situation of Example~1, the resulting augmented system is not very stiff.

For the numerical solution of this problem we use the code {\sc Radar5}.
We consider a fixed tolerance ${\it Tol} = 10^{-8}$, and we apply
the code with ${\it Atol} = {\it Rtol} = {\it Tol}$. 
We include in the mesh the first $10$ breaking points, 
\[
\tau, \ \beta, \ 2 \tau, \ \tau+\beta, \ 2 \beta, \ 3 \tau, \ 2 \tau + \beta, \ 2 \beta + \tau, \ 3 \beta, \ 4 \tau,
\]
and let the code possibly compute further breaking points.

As initial step size we choose
$h= {\it Atol}$. Similar to the first experiment for Example~1, we report the relative 
error in Table \ref{tab:1} at the final point.
We get $120$ steps (with no rejections), $854$ function evaluations, $72$ Jacobian evaluations,
$99$ LU decomposition of the matrix in the linear system and $244$ solutions of triangular systems.

The result is very similar to that of Table~\ref{tab:21}. For a fixed tolerance ${\it Tol}$,
the error is proportional to $\eps$ as long as $\eps > {\it Tol}$. It stagnates at
${\it Tol}$ for values of $\eps$ that are smaller than ${\it Tol}$.

\subsection{Example 3: a distributed delay model of chemotherapy-induced myelosuppression}
\label{sect:chemo}

We consider a model system that is proposed and studied in \cite[p.\,56]{Krzyzanski19}.
It uses ideas of the publication \cite{Humphries18}, where transit compartments in
a semi-physiological model by Friberg are replaced by a convolution integral with
a Gamma-distribution (see also \cite{Krzyzanski18}). The equations of the model are
\begin{equation}\label{example3eh}
\begin{array}{rcl}
\dot{y} (t) & = & \Bigl(\kappa  \Bigl( \displaystyle \frac{ w_0}{w(t)} \Bigr)^\gamma 
- k_s C(t)  - \kappa \Bigr) y(t)   
\\[3.5mm]
\dot{w} (t) & = & - \kappa w(t) + \kappa I \bigl( t,y(t)\bigr) 
\\[1.5mm]
\dot{A}(t) & = & \displaystyle -\frac{V_{\max} A(t)}{K_m + C(t)}, \qquad
C(t) = {\displaystyle \frac{A(t)}{V}} ,
\end{array}
\end{equation}
where $y(t)$ gives the proliferating precursor cells for granulocytes,
$w(t)$ the circulating granulocytes, and $A(t)$ is the amount of drug in the plasma.
The distributed delay term,
with a Gamma-distribution kernel, is given by
\begin{equation}\label{distridelay2eh}
I (y)(t) = \int_0^t k(t-s) y(s) \,\d s , \qquad k(t) = 
\frac {\kappa^{1-\alpha}}{\Gamma (1-\alpha)}\, t^{-\alpha} \e^{-\kappa t} .
\end{equation}
Table~\ref{tab:31eh} gives two sets of parameters, which are taken from
\cite[Table 2]{Krzyzanski19}. Concerning the units of the parameters we refer to
\cite{Krzyzanski19}. Time is measured in hours.
\begin{table}[ht] 
\begin{center}
{\small
\begin{tabular}{|c|c|c|c|c|c|c|c|}
\hline
 $\nu = 1- \alpha$ & $\kappa$ & $w_0$ & $\gamma$ & $k_s$ & $V_{\max} $ & $K_m$ & $V$\\
\hline
 \rule[0mm]{0cm}{3.0mm}%
0.964 & $\nu/47.5$ & 14.4 & 0.664 & 0.0328 & 77.2 & 16.9 & 1.35 \\
1.46 & $\nu/55.6$ & 14.4 & 0.507 & 0.0213 & 100 & 22 & 1.03 \\
  \hline
\end{tabular}}
\vskip 1.5mm
\caption{Values of the parameters for the system \eqref{example3eh}}
\label{tab:31eh} 
\end{center}
\end{table} 

We note that the third equation of \eqref{example3eh} is a scalar differential
equation for $A(t)$, which is independent of the other two variables of the system.
It can be solved analytically by separation of variables. In this way we find that
$A(t)$ is solution of the nonlinear equation
\begin{equation*}
\frac{K_m}{V_{\max} } \Bigl( \ln A(t) - \ln A_0 \Bigr)  + \frac{A(t) - A_0}{V_{\max} V} = 
- (t  -t_0 ) .
\end{equation*}
To get a better conditioning for its numerical treatment we write it as
\begin{equation}\label{equationA}
A(t) = A_0 \exp\Bigl( - \frac 1{K_mV} \bigl(A(t) - A_0 \bigr) 
- \frac{V_{\max}}{K_m}(t  -t_0 ) \Bigr) .
\end{equation}
Without changing the solution of the problem
we can replace in \eqref{example3eh} the differential equation
for $A(t)$ by the algebraic equation~\eqref{equationA}. This then leads to a
differential-algebraic system of index 1, which can be solved efficiently by {\sc Radau5}
and by {\sc Radar5}.
The formulation as a differential-algebraic system has the advantage
that during the numerical integration there is no accumulation of errors in the variable $A(t)$.
In our experience the DAE formulation is slightly more efficient
than the ODE formulation (see Table~\ref{tab:ode-dae}).

As proposed in Section~\ref{sect:tol} we add a new variable $y_{3+1}(t)$ to the
equation \eqref{example3eh}, and we apply the code {\sc Radau5} to the augmented system.
For various values of $\eps$ we use
${\it Atol} = {\it Rtol} = \eps$ for the three components of \eqref{example3eh},
${\it Atol} = {\it Rtol} = 10^{-2}\eps$ for the $y_{3+1}$-component corresponding to the
integral delay,
and ${\it Atol} = {\it Rtol} = 10^2 \eps$ for the $z$-components.
As initial step size we take $\max (\eps , 10^{-5})$.

\medskip

\begin{figure}[t]
\vspace{-1mm}
\centering
 \begin{picture}(0,0)
  \epsfig{file=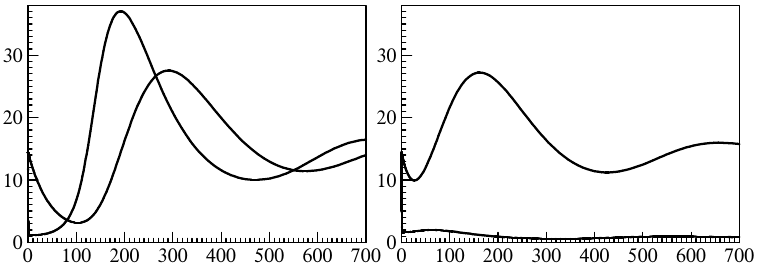}
 \end{picture}%
\begin{picture}(364.1,127.9)(  -5.8, 472.3)
  \GGGput[   50,   50](159.11,492.64){$t$}
  \GGGput[   50,   50]( 29.27,573.46){$y(t)$}
  \GGGput[   50,   50](110.45,552.50){$w(t)$}
  \GGGput[   50,   50](133.58,582.50){problem 1}
  \GGGput[   50,   50](338.34,492.64){$t$}
  \GGGput[   50,   50](210.91,498.66){$y(t)$}
  \GGGput[   50,   50](264.15,543.59){$w(t)$}
  \GGGput[   50,   50](312.81,582.50){problem 2}
 \end{picture}
\vspace{-1mm}
\caption{Solution components $y(t)$ and $w(t)$ of the system \eqref{example3eh} for the
parameter values of Table~\ref{tab:31eh}: upper row (problem 1) and lower row (problem 2).
\label{fig:solution}}
\end{figure}

\noindent{\it
First experiment.}
We consider the problem \eqref{example3eh} with parameters taken from the upper
row of Table~\ref{tab:31eh}. Initial values are
\begin{equation}\label{initialval}
y(0) = w_0,\qquad  w(0)=w_0, \qquad A(0) = A_0 = 127 .
\end{equation}
The solution components $y(t)$ and $w(t)$ are plotted in the left
picture of Figure~\ref{fig:solution}. We do not include the graph of $A(t)$, because it
monotonically decreases, rapidly approaches zero, and stays there for all $t$. 
We observe that the function $y(t)$ has an initial layer. It
rapidly decreases from $y(0)=14.4$ to a small value and then behaves smoothly.

The challenge of this problem is the fact that the value of $\alpha =0.036$ is
very small. According to Algorithm~\ref{alg_paramG} this then leads to a very small
$x_*$, and consequently to a large negative value of $M$, so that the dimension
of the augmented system becomes large. The values of $M$, $N$ (as well as $h$)
are computed by Algorithm~\ref{alg_paramG} of Section~\ref{sect:gamma-dist}.
They are presented for different values of $\eps$ in Table~\ref{tab:32eh}.

We apply the code {\sc Radau5} for
values of $\eps$ ranging between $10^{-3}$ and $10^{-10}$ on the interval
$[0,100]$. Since the Jacobian of the vector field is not banded, we apply the code
in the standard way, where the linear algebra considers the Jacobian as a full matrix.
The cpu time (measured in seconds)
is given in the row of Table~\ref{tab:32eh} marked ``cpu (full)''.
We also have implemented the treatment of the arising linear systems according to
the algorithm of Section~\ref{sect:lin-syst}, and the cpu
times are listed in the row of Table~\ref{tab:32eh} marked ``cpu (sumexp)''.
The improvement is amazing. Already for
$\eps = 10^{-3}$ the cpu time is decreased by a factor of $100$ and for
$\eps = 10^{-10}$, where the dimension of the augmented system is very large,
the cpu time is even decreased by a factor of $10^4$.
%

\medskip
\noindent{\it
Second experiment.}
We next consider
the set of parameters given in the lower row of Table~\ref{tab:31eh}.
The solution is plotted in the right
picture of Figure~\ref{fig:solution}.
Here, the value of $\alpha =-0.46$ is
negative, and the trick explained in Section~\ref{sect:gamma-dist} has to be applied.
This means that the factor $t^{-\alpha}$ in \eqref{distridelay2eh} is written as
$t^{-\alpha} = t \cdot t^{-(\alpha +1)}$ with $\alpha +1 = 0.54 \in (0,1)$, and only
the term $t^{-(\alpha +1)}$ is replaced by a sum of exponentials. This yields an
approximation of the form \eqref{sum-exp-pol} with $m_i = 1$.

\begin{table}[t] 
\begin{center}
{\small
\begin{tabular}{|c||c|c|c|c|c|c|}
\hline
 \rule[0mm]{0cm}{3.6mm}%
 $ \eps $ & $10^{-3}$ & $10^{-4}$ 
 & $10^{-6}$ & $10^{-7}$ &
 
 $10^{-9}$ & $10^{-10}$\\
\hline
 \rule[0mm]{0cm}{3.0mm}%
$M$ & $-157$ & $-268$  & $-582$ & $-783$ & $-1276$ & $-1567$ \\
$N$ & $4$ & $8$ & $20$ & $27$  & $45$ & $56$ \\
cpu (full) & $1.9\cdot 10^{-1}$ & $1.3\cdot 10^{0}$ &  $2.6\cdot 10^{1}$ &
$8.3\cdot 10^{1}$ & $5.1\cdot 10^{2}$ & $1.2\cdot 10^{3}$ \\
cpu (sumexp) & $9.6\cdot 10^{-4}$ & $2.1\cdot 10^{-3}$ & $8.9\cdot 10^{-3}$ &
$1.6\cdot 10^{-2}$ & $5.3\cdot 10^{-2}$ & $9.5\cdot 10^{-2}$ \\
  \hline
\end{tabular}}
\vskip 1.5mm
\caption{Comparison of different linear algebra solvers for the system \eqref{example3eh}
with parameters of the upper row of Table~\ref{tab:31eh}.}
\label{tab:32eh} 
\end{center}
\end{table} 

As before we apply the code {\sc Radau5} with several different values of $\eps$,
and initial values \eqref{initialval}.
We only use the option (sumexp). This time we study the error of the numerical approximation.
A reference solution is computed with a very high accuracy requirement. The exact solution
of $A(100)$ is far below round-off. The row, indicated as ``err'' in Table~\ref{tab:33eh}
shows the relative error of $y(t)$ and $w(t)$ at time $t=100$.
The values of $M$ and $N$ are computed by Algorithm~\ref{alg_paramG}. They are of
moderate size, even for very small $\eps$. We realize that the error is nicely
proportional to the accuracy requirement ${\it Tol}$.
\begin{table}[ht] 
\begin{center}
{\small
\begin{tabular}{|c||c|c|c|c|c|c|c|c|}
\hline
 \rule[0mm]{0cm}{3.6mm}%
 $ \eps $ & $10^{-3}$ & $10^{-5}$ & $10^{-7}$ & $10^{-9}$ & $10^{-11}$\\
\hline
 \rule[0mm]{0cm}{3.0mm}%
$M$ & $-17$ & $-38$ & $-67$ & $-105$ & $-150$  \\
$N$ & $13$ & $35$ & $67$ & $108$ & $156$\\
cpu (sumexp) & $4.7\cdot 10^{-4}$ & $1.7\cdot 10^{-3}$ & $5.1\cdot 10^{-3}$ & $1.5\cdot 10^{-2}$ &
$3.9\cdot 10^{-2}$\\
err & $5.3\cdot 10^{-4}$ & $1.4\cdot 10^{-5}$ & $1.5\cdot 10^{-7}$ & $3.3\cdot 10^{-9}$ &
$1.0\cdot 10^{-10}$  \\
  \hline
\end{tabular}}
\vskip 1.5mm
\caption{Study of accuracy of {\sc Radau5} for the system \eqref{example3eh}
with parameters of the lower row of Table~\ref{tab:31eh}.}
\label{tab:33eh} 
\end{center}
\end{table}

\medskip
\noindent{\it
Third experiment.}
The equations \eqref{example3eh} are an ODE formulation of the problem
of Section~\ref{sect:chemo}. When the differential equation for $A(t)$ is replaced by
the algebraic equation \eqref{equationA} we get a DAE formulation with
two differential equations and one algebraic relation.
The code {\sc Radau5} can be applied to each of these formulations, and it
is interesting to know which one performs better.
In Table~\ref{tab:ode-dae} we present for various $\eps$ the error at the final point
$t_f = 100$, as well as the required number of steps and function evaluations.
We observe that, for the
DAE formulation (the final three columns), not only the
error is smaller (excepting $\eps = 10^{-3}$), but also
the number of steps and function evaluations (and consequently the cpu time) are smaller.

\begin{table}[h] 
\begin{center}
{\small
\begin{tabular}{|c||c|c|c||c|c|c||c|c|c|} \hline
 $\eps$ & $\rule{0mm}{3mm}h$ & $M$ & $N$ & ${\rm err}_1$ & $\#{\rm st}_1$ & $\# {\rm fe}_1$ & 
 ${\rm err}_2$ & $\#{\rm st}_2$ & $\# {\rm fe}_2$ \\
\hline
\hline
\rule[0mm]{0cm}{3.5mm}%
$10^{-3}$  &  $  1.04$  &  $ -17$  &  $  13$  & 
$5.34 \cdot 10^{-4}$   &  $    25$  &  $   161$  & 
$9.54 \cdot 10^{-3}$  &  $    23$  &  $   154$   \\
$10^{-5}$  &  $  0.69$  &  $ -38$  &  $  35$  &  
$1.37 \cdot 10^{-5}$  &  $    47$  &  $    287$  & 
$9.06 \cdot 10^{-6}$  &  $    38$  &  $   262$   \\
$10^{-7}$  &  $  0.52$  &  $ -67$  &  $  67$  &  
$1.52 \cdot 10^{-7}$   &  $    80$  &  $    507$  & 
$5.47 \cdot 10^{-8}$  &  $    68$  &  $    483$   \\
$10^{-9}$  &  $  0.42$  &   $-105$  &  $ 108$  &  
$3.33 \cdot 10^{-9}$  &  $    152$  &  $   985$ &  
$4.98 \cdot 10^{-10}$  &  $   126$  &  $   945$   \\
$10^{-11}\!\!$  &  $  0.35$  &   $-150$  &  $ 156$  &  
$1.05 \cdot 10^{-10}\!\!$  &  $   309$  &  $ \!2036\!$  &
$2.41 \cdot 10^{-11}\!\!$  &  $   256$  &  $\!1933\!$   \\
\hline
\end{tabular}
}
\vskip 1.5mm
\caption{Performance of {\sc Radau5} for the problem \eqref{example3eh}
(coefficients from the lower row of Table~\ref{tab:31eh})
with $t_f=100$ and $\eps$ ranging from $10^{-3}$ to $10^{-11}$.
Here, ${\rm err}_1$, $\# {\rm st}_1$, $\# {\rm fe}_1$ indicate the relative error,
the number of steps and the number of function evaluations for an implementation
as ODE,
while ${\rm err}_2$, $\# {\rm st}_2$, $\# {\rm fe}_2$ are the numbers corresponding
to an implementation as DAE.}
\label{tab:ode-dae} 
\end{center}
\vspace{-2mm}
\end{table}

\subsection{Implementation of the approach of Section~\ref{sect:lin-syst} for solving
the linear system}
\label{sect:implement}

The code {\sc Radau5} \cite{hairer96sod} (and similarly the code
{\sc Radar5} \cite{guglielmi01irm}) consists of three files: {\sc Radau5.f} 
({\sc Radar5.f}) contains the time integrator, 
{\sc decsol.f} the linear algebra subroutines, and
{\sc dc\_decsol.f} ({\sc dc\_decdel.f}) contains the subroutines that link the time integrator with the linear algebra subroutines. When using the linear algebra approach of
Section~\ref{sect:lin-syst} one only has to replace the file {\sc dc\_decsol.f}
(or {\sc dc\_decdel.f}) by {\sc dc\_sumexp.f}
(or {\sc dc\_sumexpdel.f})\footnote{The subroutines {\sc dc\_decdel.f} and {\sc dc\_sumexpdel.f}
are restricted to the case where $M$ is a (possibly singular) diagonal matrix.}. 
The other files 
need not be touched.
In the driver the Jacobian has to be defined as follows:
first of all, if $d$ is the dimension of $y$ in the original system and $n$ the number 
of summands in \eqref{sum-exp-pol}, we have to define the dimension of the system as
\[
N = \biggl\{ \begin{array}{ll}d + 1 + n + 2 & \hbox{if all $m_i=0$ in \eqref{sum-exp-pol}} \\[1.1mm]
d + 1 + 2n + 2 & \hbox{if all $m_i=1$ in \eqref{sum-exp-pol}. }
\end{array} 
\]
The first $d$ components correspond to those of the original system \eqref{problem},
the $(d+1)$th component is for the integral delay \eqref{newvariable}, and the next $n$ (or $2n$)
components are for the variables $z_{i,j}$. The final $2$ components are not used in the
right-hand side function of the problem. They are used only for the Jacobian.

The Jacobian has to be declared as banded with upper bandwidth $\mbox{\sc mujac}=d$
and lower bandwith $\mbox{\sc mljac}=\max (2-d,0)$. The array {\sc fjac} then has
$N$ columns and $\max (d+1,3)$ rows.
In the left upper $(d+1)\times (d+1)$ matrix, the derivative of
the vector field $f$ (augmented by \eqref{newvariable}) with
respect to $(y,y_{d+1})$ has to be stored, 
in column $d+2$ the derivative
of $f$ with respect to the integral term $I$, and in column $d+3$ the gradient of
$g(t,y)$ with respect to $y$. For the case that all $m_i=0$, the coefficients
$\{c_{i0}\, |\, i=1,\ldots n \}$ are stored in the first row starting at position $d+4$,
the exponents $\{\gamma_i\, |\,  i=1,\ldots n \}$ in the second row,
and zeros in the third row. Finally, for the case that all $m_i=1$, the coefficients
$\{(c_{i0}, c_{i1})\, |\, i=1,\ldots n \}$ are stored in the first row starting at
position $d+4$, the exponents $\{(\gamma_i , \gamma_i)\, |\,  i=1,\ldots n \}$ in the second row,
and the subdiagonal $(1,0,1,0, \ldots , 1,0)$ of \eqref{jacobian} in the third row.

The codes {\sc Radau5} and {\sc Radar5} together with the files
{\sc decsol.f}, {\sc dc\_decsol.f}, {\sc dc\_sumexp.f}
together with drivers for the examples of this article are available at the address 
\url{http://www.unige.ch/~hairer/software.html}.

\section*{Conclusions}

In this article we have described an efficient methodology to deal with important classes of distributed
delays, by using {\sc Radau5} or {\sc Radar5}, which are codes for the numerical integration
of a large class of stiff and differential-algebraic
equations, and with discrete delays (for {\sc Radar5}). 

The methodology is based on a suitable expansion of the kernels in terms
of exponential functions, which allows to transform the distributed delay into a set of ODEs or DDEs.
Our numerical experiments confirm the efficiency of this methodology which
allows to get
a much more versatile code, able to deal with both discrete and distributed memory effects,
a feature that
makes it very appealing in disciplines where these kind of delays naturally arise,
as pharmacodynamics
and pharmacokinetics.

\subsection*{Acknowledgments}

Nicola Guglielmi acknowledges that his research was supported by funds from the Italian 
MUR (Ministero dell'Universit\`a e della Ricerca) within the 
`
PRIN 2022 Project ``Advanced numerical methods for time dependent parametric partial differential equations with applications'' and the   Pro3 joint project entitled
``Calcolo scientifico per le scienze naturali, sociali e applicazioni: sviluppo metodologico e tecnologico''.
He is also affiliated to the INdAM-GNCS (Gruppo Nazionale di Calcolo Scientifico).

Ernst Hairer acknowledges the support of the Swiss National Science Foundation, grant No.200020 192129.


\end{document}